\magnification=\magstep1
\hsize=16.5 true cm 
\vsize=25 true cm

\font\bfff=cmbx10 scaled \magstep2

\font\bffgg=cmbx10 scaled \magstep4
\font\smc=cmcsc10 
\parindent0cm
           
\def\Bbb#1{\hbox{\boldmas #1}} 
 
\expandafter\edef\csname amssym.def
\endcsname{%
       \catcode`\noexpand\@=\the\catcode`\@\space}

\catcode`\@=11

\def\undefine#1{\let#1\undefined}
\def\newsymbol#1#2#3#4#5{\let\next@\relax
 \ifnum#2=\@ne\let\next@\msafam@\else
 \ifnum#2=\tw@\let\next@\msbfam@\fi\fi
 \mathchardef#1="#3\next@#4#5}
\def\mathhexbox@#1#2#3{\relax
 \ifmmode\mathpalette{}{\m@th\mathchar"#1#2#3}%
 \else\leavevmode\hbox{$\m@th\mathchar"#1#2#3$}\fi}
\def\hexnumber@#1{\ifcase#1 0\or 1\or 2\or 3\or 4\or 5\or 6\or 7\or 8\or
 9\or A\or B\or C\or D\or E\or F\fi}

\font\tenmsa=msam10
\font\sevenmsa=msam7
\font\fivemsa=msam5
\newfam\msafam
\textfont\msafam=\tenmsa
\scriptfont\msafam=\sevenmsa
\scriptscriptfont\msafam=\fivemsa
\edef\msafam@{\hexnumber@\msafam}
\mathchardef\dabar@"0\msafam@39
\def\dashrightarrow{\mathrel{\dabar@\dabar@\mathchar"0\msafam@4B}}
\def\dashleftarrow{\mathrel{\mathchar"0\msafam@4C\dabar@\dabar@}}

\def\ulcorner{\delimiter"4\msafam@70\msafam@70 }
\def\urcorner{\delimiter"5\msafam@71\msafam@71 }
\def\llcorner{\delimiter"4\msafam@78\msafam@78 }
\def\lrcorner{\delimiter"5\msafam@79\msafam@79 }
\def\yen{{\mathhexbox@\msafam@55}}
\def\checkmark{{\mathhexbox@\msafam@58}}
\def\circledR{{\mathhexbox@\msafam@72}}
\def\maltese{{\mathhexbox@\msafam@7A}}

\font\tenmsb=msbm10
\font\sevenmsb=msbm7
\font\fivemsb=msbm5
\newfam\msbfam
\textfont\msbfam=\tenmsb
\scriptfont\msbfam=\sevenmsb
\scriptscriptfont\msbfam=\fivemsb
\edef\msbfam@{\hexnumber@\msbfam}
\def\Bbb#1{{\fam\msbfam\relax#1}}
\def\widehat#1{\setbox\z@\hbox{$\m@th#1$}%
 \ifdim\wd\z@>\tw@ em\mathaccent"0\msbfam@5B{#1}%
 \else\mathaccent"0362{#1}\fi}

\def\widetilde#1{\setbox\z@\hbox{$\m@th#1$}%
 \ifdim\wd\z@>\tw@ em\mathaccent"0\msbfam@5D{#1}%
 \else\mathaccent"0365{#1}\fi}
\font\teneufm=eufm10
\font\seveneufm=eufm7
\font\fiveeufm=eufm5
\newfam\eufmfam
\textfont\eufmfam=\teneufm
\scriptfont\eufmfam=\seveneufm
\scriptscriptfont\eufmfam=\fiveeufm

\newsymbol\risingdotseq 133A
\newsymbol\fallingdotseq 133B
\newsymbol\complement 107B
\newsymbol\nmid 232D
\newsymbol\rtimes 226F
\newsymbol\thicksim 2373

\font\eightmsb=msbm8   \font\sixmsb=msbm6   \font\fivemsb=msbm5
\font\eighteufm=eufm8  \font\sixeufm=eufm6  \font\fiveeufm=eufm5
\font\eightrm=cmr8     \font\sixrm=cmr6     \font\fiverm=cmr5
\font\eightbf=cmbx8    \font\sixbf=cmbx6    
      \font\eighti=cmmi8   \font\sixi=cmmi6
\font\ninesy=cmsy9     \font\eightsy=cmsy8  \font\sixsy=cmsy6
     \font\eightit=cmti8  
     \font\eightsl=cmsl8  
     \font\eighttt=cmtt8

\font\eightsmc=cmcsc8
\newskip\ttglue
\newfam\smcfam
\def\eightpoint{\def\rm{\fam0\eightrm}%
  \textfont0=\eightrm \scriptfont0=\sixrm \scriptscriptfont0=\fiverm
  \textfont1=\eighti \scriptfont1=\sixi \scriptscriptfont1=\fivei
  \textfont2=\eightsy \scriptfont2=\sixsy \scriptscriptfont2=\fivesy
  \textfont3=\tenex \scriptfont3=\tenex \scriptscriptfont3=\tenex
  \def\smc{\fam\smcfam\eightsmc}
  \textfont\smcfam=\eightsmc          
\textfont\eufmfam=\eighteufm              \scriptfont\eufmfam=\sixeufm
     \scriptscriptfont\eufmfam=\fiveeufm
\textfont\msbfam=\eightmsb            \scriptfont\msbfam=\sixmsb
     \scriptscriptfont\msbfam=\fivemsb
\def\it{\fam\itfam\eightit}%
  \textfont\itfam=\eightit
  \def\sl{\fam\slfam\eightsl}%
  \textfont\slfam=\eightsl
  \def\bf{\fam\bffam\eightbf}%
  \textfont\bffam=\eightbf \scriptfont\bffam=\sixbf
   \scriptscriptfont\bffam=\fivebf
  \def\tt{\fam\ttfam\eighttt}%
  \textfont\ttfam=\eighttt
  \tt \ttglue=.5em plus.25em minus.15em
  \normalbaselineskip=9pt
  \def\MF{{\manual opqr}\-{\manual stuq}}%
  \let\big=\eightbig
  \setbox\strutbox=\hbox{\vrule height7pt depth2pt width\z@}%
  \normalbaselines\rm}
\def\eightbig#1{{\hbox{$\textfont0=\ninerm\textfont2=\ninesy
  \left#1\vbox to6.5pt{}\right.\n@space$}}}

\catcode`@=13 
\centerline{\bffgg On decompositions of the real line}
\bigskip\medskip
\centerline{\bfff Gerald Kuba}
\bigskip\bigskip\medskip
\vbox{\eightpoint {\bf Abstract.}  
Let $\,X_t\,$ be a totally disconnected subset of $\,{\Bbb R}\,$
for each $\,t\in{\Bbb R}\,$.
We construct a partition $\;\{\,Y_t\;|\;t\in{\Bbb R}\,\}\;$
of $\,{\Bbb R}\,$ into nowhere dense Lebesgue null sets $\,Y_t\,$
such that for every $\,t\in{\Bbb R}\,$
there exists an increasing homeomorphism
from $\,X_t\,$ onto $\,Y_t\,$. In particular, 
the real line can be partitioned into $\,2^{\aleph_0}\,$ 
Cantor sets and also into 
$\,2^{\aleph_0}\,$ mutually non-homeomorphic compact subspaces.
Furthermore we prove that  
for every cardinal number $\,\kappa\,$
with $\,2\leq \kappa\leq 2^{\aleph_0}\,$ the real line
(as well as the Baire space $\,{\Bbb R}\setminus{\Bbb Q}\,$)
can be partitioned into exactly
$\,\kappa\,$ {\it homeomorphic} Bernstein sets
and also into exactly 
$\,\kappa\,$ {\it  mutually non-homeomorphic} Bernstein sets.         
We also investigate partitions of $\,{\Bbb R}\,$
into Marczewski sets, including the 
possibility that they are 
Luzin sets or Sierpi\'nski sets.
}
\bigskip\medskip
{\bf 1.~Introduction.} 
The cardinal number (the {\it size}) of a set $\,S\,$
is denoted by $\,|S|\,$. So 
$\,|{\Bbb R}|=|{\Bbb D}|={\bf c}=2^{\aleph_0}\,$ where 
$\;{\Bbb D}\,=\,
\{\,\sum_{n=1}^\infty a_n\,3^{-n}\;\,|\,\;a_n\in\{0,2\}\,\}\;$
is the {\it Cantor ternary set}.
\medskip
We are interested in partitions of $\,{\Bbb R}\,$ from a specific topological point 
of view. Consider for $\,\emptyset\not=A\subset{\Bbb R}\,$ 
the following statement.
\smallskip
(1.1)\quad{\it There exists a partition $\,{\cal P}\,$ of $\,{\Bbb R}\,$
such that every $\,P\in{\cal P}\,$ is homeomorphic with $\,A\,$.}
\smallskip
Naturally, if $\,A\,$ is a subgroup of $\,({\Bbb R},+)\,$ then 
the quotient group $\,{\Bbb R}/A\,$ is a suitable partition $\,{\cal P}\,$.
If $\,A\,$ is compact then (1.1) is only possible 
if $\,A\,$ is totally disconnected.
Because, by applying Sierpi\'nski's theorem [2] 6.1.27, 
the only partition of $\,{\Bbb R}\,$ into {\it countably many} closed sets is
$\,\{{\Bbb R}\}\,$. However, if $\,A\,$ is a totally disconnected 
subset of $\,{\Bbb R}\,$, in particular if $\,A={\Bbb D}\,$, 
then (1.1) holds with $\,|{\cal P}|={\bf c}\,$.
Moreover, the following is true.
\medskip
{\bf Theorem 1.} {\it Assume $\,\emptyset\not=X_t\subset{\Bbb R}\,$ 
for every index $\,t\in{\Bbb R}\,$.
Then there exists a partition $\;\{\,Y_t\;|\;t\in{\Bbb R}\,\}\;$
of $\,{\Bbb R}\,$ into nowhere dense Lebesgue null sets $\,Y_t\,$
such that for every $\,t\in{\Bbb R}\,$
there is an increasing bijection $\,\varphi_t\,$
from $\,X_t\,$ onto $\,Y_t\,$ and $\,\varphi_t\,$ is a homeomorphism 
in case that $\,X_t\,$ is totally disconnected.}
\medskip                   
The restricton to subspaces $\,X_t\,$ of $\,{\Bbb R}\,$ is not very strong
since (in view of [2] 6.2.16) every zero-dimensional 
second countable space can be embedded into $\,{\Bbb D}\,$.
Therefore (considering [2] 6.2.9) and since (in view of [6] 1.2 and 
Proposition 2) there exist exactly 
$\,{\bf c}\,$ totally disconnected compact metrizable 
spaces up to homeomorphism, from Theorem 1 we derive 
\medskip
{\bf Corollary 1.} {\it The real line can be partitioned into 
$\,{\bf c}\,$ compact sets which form a complete system 
of representatives of all
totally disconnected compact metrizable spaces.}
\medskip
If in Theorem 1 all sets $\,X_t\,$ are identical and infinite
(in particular, if $\,X_t={\Bbb D}\,$) 
then the statement about the mappings $\,\varphi_t\,$
can be sharpened in the following way.
\medskip
{\bf Theorem 2.} {\it If $\,X\subset{\Bbb R}\,$ is infinite
then there is a nowhere dense Lebesgue null set $\,Y\subset{\Bbb R}\,$
and an increasing bijection $\,\varphi\,$
from $\,X\,$ onto $\,Y\,$, 
which is a homeomorphism when $\,X\,$ is totally disconnected,
such that $\,{\Bbb R}\,$ can be partitioned into $\,{\bf c}\,$
translates $\,\xi+Y\,$ of $\,Y\,$.}
\medskip
In the light of Theorems 1 and 2 
one may ask about partitions of $\,{\Bbb R}\,$ into
{\it dense} and {\it non-measurable} subsets 
and about partitions $\,{\cal P}\,$ with $\,|{\cal P}|<{\bf c}\,$.
Both question are answered simultaneously by the 
following theorem about {\it Bernstein sets.}
(As usual, $\,B\subset{\Bbb R}\,$ is a Bernstein set 
if and only if neither $\,B\,$ nor $\,{\Bbb R}\setminus B\,$ 
contains an uncountable closed set.)
\medskip
{\bf Theorem 3.} {\it If $\,2\leq \kappa\leq {\bf c}\,$ then $\,{\Bbb R}\,$
can be partitioned into exactly $\,\kappa\,$ {\it homeomorphic} Bernstein sets
and also into exactly 
$\,\kappa\,$ {\it  mutually non-homeomorphic} Bernstein sets.}
\medskip
Notice that, by a straightforward application of Easton's theorem [3] 15.18,
{\it it is consistent with {\rm ZFC} set theory that there are 
$\,{\bf c}\,$ cardinal numbers $\,\kappa\,$ with $\,\aleph_0<\kappa<{\bf c}\,$.}
\eject
\bigskip
{\bf 2.~Cantor sets.} A nonempty set
$\,C\subset{\Bbb R}\,$ is a {\it Cantor set} if and only if $\,C\,$ 
is compact and does not contain nondegenerate intervals 
or isolated points. It is well-known that $\,C\subset{\Bbb R}\,$
is a Cantor set if and only if $\,C\,$ 
homeomorphic to $\,{\Bbb D}\,$. For abbreviation, 
a {\it Cantor null set} is a Cantor set which is a 
Lebesgue null set. In particular, if $\,x,t\in{\Bbb R}\,$ and $\,t>0\,$
then $\;x+t\cdot{\Bbb D}\;$ is a Cantor null set.
Furthermore, for a Cantor set $\,C\,$ let $\,C^*\,$
denote the set of all $\,x\in {\Bbb R}\,$ 
such that $\;C\cap\,]x,x+\epsilon[\,\not=\emptyset\;$ 
and $\;C\cap\,]x-\epsilon,x[\,\not=\emptyset\;$ 
for every $\,\epsilon>0\,$.
Clearly, $\,C^*\subset C\,$ and $\,C\setminus C^*\,$ is countable.
For the proof of Theorem 1 we need three lemmas.
\medskip
{\bf Lemma 1.} {\it If $\,S\subset{\Bbb R}\,$ is totally disconnected
and $\,C\,$ is a Cantor set and $\,s\in S\,$ and $\,a\in C^*\,$ 
then there exists a strictly increasing function $\,f\,$ from 
$\,{\Bbb R}\,$ into $\,C\,$ such that $\;f(s)=a\;$ and 
$\,f\,$ restricted to the domain 
$\,S\,$ is a homeomorphism from $\,S\,$ onto the subspace $\,f(S)\,$
of $\,C\,$.}          
\medskip
{\it Proof.} Certainly we can fix a countable
dense subset $\,D\,$ of $\,{\Bbb R}\,$ disjoint from $\,S\,$.
Let $\,{\cal U}_1\,$ and $\,{\cal U}_2\,$ be the family of the components
of the open subspace $\;[\min C,a]\setminus C\;$ and 
$\;[a,\max C]\setminus C\;$ of $\,{\Bbb R}\,$, respectively.
The set $\;{\cal U}_1\cup{\cal U}_2\;$ is naturally ordered 
such that the sets $\,{\cal U}_1\cup{\cal U}_2\,$
and $\,{\cal U}_1\,$ and $\,{\cal U}_2\,$ are order-isomorphic with $\,{\Bbb Q}\,$
and hence order-isomorphic with the sets
$\,D\,$ and $\;D_1\,:=\;]{-\infty,s}]\cap D\;$ and 
$\;D_2\,:=\,[s,\infty[\cap D\,$. Thus we can define 
a bijection $\,\varphi\,$
from $\,D\,$ onto $\,{\cal U}_1\cup{\cal U}_2\,$ 
such that $\;\varphi(D_1)={\cal U}_1\;$ 
(and hence $\;\varphi(D_2)={\cal U}_2\,$)
and for $\,x,y\in D\,$ we have
$\,x<y\,$ if and only if $\;\sup\varphi(x)<\inf\varphi(y)\,$.
Finally, define a function $\,f\,$ 
from $\,{\Bbb R}\,$ into $\,C\,$ by
\smallskip
\centerline{$f(x)\;=\;\sup\,\bigcup\{\,\varphi(d)\;\,|\,\;d\in D\;\land\;d<x\,\}
\qquad(x\in{\Bbb R})\,$.}
\smallskip
Clearly, $\,f\,$ is strictly increasing with
$\,f(s)=a\,$ and
$\,f\,$ restricted to the domain $\;{\Bbb R}\setminus D\;$ is
a homeomorphism from $\;{\Bbb R}\setminus D\;$ onto $\,C^*\,$,
{\it q.e.d.}
\medskip
By similar arguments we clearly obtain the following lemma.
\medskip
{\bf Lemma 2.} {\it If $\,C_i\subset{\Bbb R}\,$ is a Cantor set 
and $\,a_i\in C_i^*\,$ for $\,i\in\{1,2\}\,$ 
then there exists an increasing homeomorphism $\,g\,$ from 
$\,{\Bbb R}\,$ onto $\,{\Bbb R}\,$ such that $\,g(a_1)=a_2\,$ and 
$\,g(C_1)=C_2\,$.}          
\medskip
{\bf Lemma 3.} {\it There exists a family $\,{\cal H}\,$ 
of pairwise disjoint Cantor null sets such that 
every nonempty open interval contains $\,{\bf c}\,$ 
sets from the family $\,{\cal H}\,$.}
\medskip                 
In order to settle Lemma 3 we first prove 
the following useful observation.
\smallskip
(2.1)\quad {\it Every Cantor set can be partitioned into $\,{\bf c}\,$
Cantor sets.}
\smallskip
To verify (2.1) consider the subset
$\;D\,:=\,
\{\,\sum_{n=1}^\infty a_n\,3^{-2n}\;\,|\,\;a_n\in\{0,2\}\,\}\;$
of $\,{\Bbb D}\,$ which is clearly a Cantor set. Obviously, 
$\;\{\,x+3\cdot D\;|\;x\in D\,\}\;$
is a partition of $\,{\Bbb D}\,$ into $\,{\bf c}\,$ Cantor sets.
Now let $\,C\,$ be a Cantor set and 
let $\,f\,$ be a homeomorphism from 
$\,{\Bbb D}\,$ onto $\,C\,$. Then $\;\{\,f(x+3\cdot D)\;|\;x\in D\,\}\;$
is a partition of $\,C\,$ into $\,{\bf c}\,$ Cantor sets, {\it q.e.d.}
\medskip
Now in order to prove Lemma 3 let 
$\;\{\,I_n\;|\;n\in{\Bbb N}\,\}\;$ be a basis of the space $\,{\Bbb R}\,$ 
consisting of nonempty open intervals. 
Clearly, if $\,u<v\,$ then for some $\,x,t\in{\Bbb R}\,$ with $\,t>0\,$
the Cantor null set $\;x+t\cdot{\Bbb D}\;$ lies in $\,]u,v[\,$.
Therefore, step by step we can choose a Cantor null set $\,C_1\subset I_1\,$
and a Cantor null set $\;C_2\,\subset\,I_2\setminus C_1\;$
and a Cantor null set $\;C_3\,\subset\,I_3\setminus (C_1\cup C_2)\;$
and so on, finally obtaining pairwise disjoint 
Cantor null sets $\;C_n\;(n\in{\Bbb N})\;$ with $\;C_n\subset I_n\;$ 
for every $\,n\in{\Bbb N}\,$. By applying (2.1) to each Cantor null set $\,C_n\,$
we obtain a family $\;{\cal H}\,$ as desired 
and this concludes the proof of Lemma 3.
\medskip\smallskip
{\bf 3.~Proof of Theorem 1.}
In the following, a {\it minimal well-ordering} of an infinite set $\,X\,$ 
is a linear ordering $\,\preceq\,$ of $\,X\,$ 
such that $\;|\{\,y\;|\;y\prec x\,\}|<|X|\;$ 
for every $\,x\in X\,$ and that 
every nonempty subset of $\,X\,$ has a $\preceq$-minimum. 
(In other words, $\,X\,$ equipped with $\,\preceq\,$ 
is order-isomorphic with the naturally ordered 
set of all ordinal numbers smaller than
the initial ordinal number $\,|X|\,$.)
Now we regard $\,{\Bbb R}\,$ both as 
the basic linearly ordered space 
in Theorem 1 and as an index set, assuming
that $\;\xi_t\in X_t\subset{\Bbb R}\;$ 
for every index $\,t\in{\Bbb R}\,$.
Fix a real number $\,\theta\in{\Bbb D}^*\,$ and consider any index $\,t\in{\Bbb R}\,$.
\smallskip
If $\,X_t\,$ is not totally disconnected then,
by applying Lemma 1 with $\,S=\{x_t\}\,$, 
we can define  a strictly increasing function 
$\;f_t:\,{\Bbb R}\to{\Bbb D}\;$ with $\,f_t(\xi_t)=\theta\,$.
If $\,X_t\,$ is totally disconnected then, 
by applying Lemma 1 with $\,S=X_t\,$,  
there is an increasing homeomorphism $\,f_t\,$ from 
$\,X_t\,$ onto a subspace of $\,{\Bbb D}\,$ with $\,f_t(\xi_t)=\theta\,$.
So in both cases, $\,f_t(X_t)\subset{\Bbb D}\,$ and
$\,f_t(\xi_t)=\theta\,$.
\smallskip
Let $\,{\cal H}\,$ be a family as in Lemma 3.
In particular, $\,|{\cal H}|={\bf c}\,$.
Let $\,\preceq\,$ be a minimal well-ordering
of $\,{\Bbb R}\,$. By induction we define triples $\;(a_t,A_t,{\cal H}_t)\;$ 
for every index $\,t\in{\Bbb R}\,$ as follows.
Let $\,y\,$ be an arbitrary real number and assume
that for each $\,x\prec y\,$ a 
triple $\;(a_x,A_x,{\cal H}_x)\;$ is already defined 
where $\,a_x\in{\Bbb R}\,$ and $\,{\cal H}_x\,$ is a countable subfamily 
of $\,{\cal H}\,$ and
$\,A_x\subset \bigcup{\cal H}_x\,$.
(Notice that we abstain from 
burdening the inductive assumption with additional conditions
that will be enforced by the inductive step anyhow.
For example, $\,A_u\cap A_v=\emptyset\,$ will be enforced 
for arbitrary $\,u\prec v\,$ without assuming it for $\,u\prec v\prec y\,$.)
\smallskip
To carry out the inductive step 
put $\;U_y\,:=\,\bigcup_{x\prec y}(\{a_x\}\cup A_x)\;$
and $\;{\cal F}_y\,:=\,\bigcup_{x\prec y}{\cal H}_x\;$ 
and let $\,{\cal G}_y\,$ be the family of all 
Cantor sets in $\,{\cal H}\,$ which intersect the set 
$\;\{\,a_x\;|\;x\prec y\,\}\,$.
Notice that $\,|{\cal F}_y|<{\bf c}\,$ by definition and $\,|{\cal G}_y|<{\bf c}\,$
since the Cantor sets in $\,{\cal H}\,$ are pairwise disjoint.
Since $\;|\{\,a_x\;|\;x\prec y\,\}|<{\bf c}\;$
and $\;\bigcup_{x\prec y}A_x\,\subset\,\bigcup{\cal F}_y\;$ and
$\;|{\cal H}\setminus {\cal F}_y|={\bf c}\,$,
the set $\;{\Bbb R}\setminus U_y\;$ is not empty 
and hence we can define $\;a_y:=\min_\preceq({\Bbb R}\setminus U_y)\,$.
In order to continue the definition of the triple
$\;(a_y,A_y,{\cal H}_y)\;$ we apply Lemma 3 
and choose Cantor null sets $\,R_n,L_n\,$ 
in the family $\;{\cal H}\setminus({\cal F}_y\cup{\cal G}_y)\;$
such that $\;R_n\subset\,]a_y+2^{-n-1},a_y+2^{-n}[\;$
and $\;L_n\subset\,]a_y-2^{-n},a_y-2^{-n-1}[\;$
for every $\,n\in{\Bbb N}\,$. Then we put 
\smallskip
\centerline{$\;{\cal H}_y\;:=\;\{\,L_n\;|\;n\in{\Bbb N}\,\}\cup\{\,R_n\;|\;n\in{\Bbb N}\,\}\,$.}
\smallskip
It is evident that $\;C_y\,:=\,\{a_y\}\cup\bigcup{\cal H}_y\;$ 
is compact and totally disconnected and dense in itself.
And of course $\;\bigcup{\cal H}_y\;$ is a Lebesgue null set. 
Thus $\,C_y\,$ is a Cantor null set 
and $\,a_y\in C_y^*\,$.
By applying Lemma 2 we define a homeomorphism $\,g_y\,$ from 
$\,{\Bbb D}\,$ onto $\,C_y\,$ with $\,g_y(\theta)=a_y\,$ or, equivalently, 
$\,g_y(f_y(\xi_y))=a_y\,$. Then, concluding 
the definition of the triple $\;(a_y,A_y,{\cal H}_y)\;$
and finishing the inductive step,
we put $\;A_y\,:=\,g_y(f_y(X_y\setminus\{\xi_t\}))\,$.
\smallskip
In this way we obtain for each index $\,t\in{\Bbb R}\,$ 
a set $\;Y_t\,:=\,\{a_t\}\cup A_t\;$ 
together with a homeomorphism $\,g_t\,$ from $\,{\Bbb D}\,$ onto some Cantor null
set such that $\,g_t(f_t(X_t))=Y_t\,$.
The sets $\,Y_t\,$ are pairwise disjoint and,
having used $\,\min_\preceq(\cdot)\,$ as a choice function 
for defining the points $\,a_t\,$, the sets $\,Y_t\,$ cover 
the whole set $\,{\Bbb R}\,$. Thus by finally defining the increasing 
bijection $\,\varphi_t\,$ from $\,X_t\,$ onto $\,Y_t\,$ via
$\;\varphi_t\;:=\,g_t\circ f_t\;$ the proof of Theorem 1 is finished.
\medskip\smallskip
{\bf 4.~Proof of Theorem 2.} Let $\,C\subset{\Bbb R}\,$ be a Cantor set 
{\it algebraically independent} over $\,{\Bbb Q}\,$.
(Such a set $\,C\,$ exists by either applying Mycielski's theorem [8] 6.5
or by easily selecting $\,C\,$ from the {\it von Neumann numbers}
as defined in [8].) Naturally (cf.~[4] Ch.~8, Theorem 2) 
$\,C\,$ must be a Lebesgue null set. 
Apply Lemma 2 to define 
an increasing homeomorphism $\;g:\,{\Bbb R}\to{\Bbb R}\;$ with $\,g({\Bbb D})=C\,$.
With $\,X=X_0\,$ define $\,f_t\,$ for $\,t=0\,$ 
as in the previous proof. Then $\;\varphi\,:=\,g\circ f_0\;$
and $\;Y:=g(f_0(X))\;$ are as depicted in Theorem~2
and hence the proof of Theorem 2 is concluded by 
finding a set $\,S\subset{\Bbb R}\,$ such that 
$\;\{\,\xi+Y\;|\;\xi\in S\,\}\;$ is a partition 
of $\,{\Bbb R}\,$ of size $\,{\bf c}\,$. This can be accomplished  
by applying the following theorem which is interesting of its own
and will be needed in Section 9.
(In accordance with [7], $\,X\,$ is an {\it independent} subset
of an abelian group if and only if $\,0\not\in X\,$ and
$\;k_i x_i=0\;(1\leq i\leq n)\;$ follows from 
$\;\sum_{i=1}^n k_i x_i\,=\,0\;$ whenever 
$\,k_i\in {\Bbb Z}\,$ and $\,\{x_1,...,x_n\}\subset X\,$ is of size $\,n\in{\Bbb N}\,$.)
\medskip\smallskip
{\bf Theorem 4.} {\it Let $\,G\,$ be an abelian group and 
$\,I\,$ an index set of size $\,|G|\,$
and let $\,Y\,$ be an infinite independent 
subset of $\,G\,$ without elements of order $\,2\,$.
Assume $\,\emptyset\not=X_i\subset Y\,$ for every $\,i\in I\,$
where $\;|\{\,i\in I\;|\;|X_i|=|Y|\,\}|=|I|\,$.
Then there is a mapping $\,\rho\,$ from $\,I\,$ 
into $\,G\,$ such that $\;G=\,\bigcup_{i\in I}(\rho(i)\!+\!X_i)\;$
and $\;(\rho(i)\!+\!X_i)\cap(\rho(j)\!+\!X_j)=\emptyset\;$ 
whenever $\,i\not=j\,$.}
\medskip\smallskip
{\it Remark.} By applying Theorem 4, in Theorem 2 
we can choose $\,Y:=X\,$ in case that $\,X\,$ 
is nowhere dense and measurable and 
{\it linearly independent} in the vector space $\,{\Bbb R}\,$  over $\,{\Bbb Q}\,$.
\vfill\eject
It is also worth mentioning that a special case in Theorem~1 
can alternatively be settled by applying Theorem 4.
Indeed, by using 
an algebraically independent Cantor set,
from Lemma 1 and
Theorem 4 we derive Theorem 1 under the (strong) restriction 
that either $\;|X_t|={\bf c}\;$ for $\,{\bf c}\,$ indices $\,t\in{\Bbb R}\,$,
or all sets $\,X_t\,$ are identical and infinite.
But even if we ignore this restriction,  
the proof of Theorem 1 still has the benefit 
of being purely topological, whereas Theorem 4 is purely algebraic. 
For example, a natural and plain adaption of the proof of Theorem 1
(and of Lemma 3)
shows that {\it if $\,\{\,X_t\;|\;t\in{\Bbb R}\,\}\;$ is any collection of 
zero-dimensional second countable spaces then 
every dense-in-itself Polish space
has a partition $\,\{\,Y_t\;|\;t\in{\Bbb R}\,\}\;$
where $\,Y_t\,$ is homeomorphic with $\,X_t\,$ for every $\,t\in{\Bbb R}\,$.}
\medskip\smallskip
{\bf 5.~Proof of Theorem 4.} First of all we point out that 
for all subsets $\,A,B\,$ of the 
independent set $\,Y\subset G\,$ (not containing 
elements of order 2) the following is true.
\smallskip
(5.1)\quad{\it $\,|\,\{\,t\in G\;\,|\,\;\theta\in t\!+\!A\;\land
(t\!+\!A)\cap B\not=\emptyset\,\}\,|\,\leq \,3\;\;$ 
for every $\;\theta\,\in\,G\setminus B\,$.}
\smallskip
(5.2)\quad{\it $\,|\,\{\,t\!+\!A\;\,|\,\;t\in G\;\land\;\theta\in t\!+\!A\,\}|
\,=\,|A|\;\;$ 
for every $\;\theta\in G\,$.}
\smallskip
To verify (5.1) assume $\,0\not=t\in G\,$. Then 
$\,\theta\in t\!+\!A\,$ 
if and only if $\,\theta=t+\alpha\,$ for some 
$\;\alpha\,\in\,A\setminus\{\theta\}\,$.
Furthermore, $\;(t+A)\cap B\not=\emptyset\;$ if and only if 
$\;t+a=b\;$ for some $\,a\in A\,$ and $\,b\in B\,$ with $\,a\not=b\,$.
Therefore (5.1) holds because (in view of $\,y+y\not=0\,$ for every 
$\,y\in Y\,$) for
$\;{\Bbb D}elta_\theta\,:=\,\{\,(\alpha,a,b)\in A\!\times\! A\! \times\! B\;\,|\,\;
\alpha\not=\theta=\alpha-a+b\,,\;a\not=b\,\}\;$ 
it is plain that 
$\,|{\Bbb D}elta_\theta|\leq 2\,$ whenever $\,\theta\not\in B\,$. 
Furthermore, while $\,|\;\{\,t\in G\;|\;\theta\in t\!+\!A\,\}|=|A|\;$
is trivially true, (5.2) is trivial for $\,|A|\leq 1\,$ and
can easily be verified for $\,|A|\geq 2\,$. (For
distinct $\,a_1,a_2\in A\,$ we derive
a contradiction with the independence of $\,A\,$
from the assumption $\;(\theta-a_1)+A=(\theta-a_2)+A\;$
since $\,a+a\not=0\,$ for every $\,a\in A\,$.)
\medskip
In order to prove Theorem 4 it is enough to treat the case $\,|Y|=|G|\,$.
Because if $\,|Y|<|G|\,$ then we consider the 
subgroup $\,H\,$ generated by $\,Y\,$ and let $\,H\,$ play the 
leading role. After having settled Theorem 4 for $\,H\,$ instead of $\,G\,$ 
with $\,I=H\,$, we consider the quotient group $\,G/H\,$
where $\,|G/H|=|G|\,$ (since $\,|H|=|Y|<|G|\,$).
Then we define a new index set $\,I'\,$ by 
$\;I'\,=\,(G/H)\times H\;$ and consider that 
$\,G/H\,$ consists of $\,|G|\,$ translates of $\,H\,$. 
\medskip
Now assume $\,|Y|=|G|=:\kappa\;$ and
regard the index set $\,I=G\,$
equipped with a minimal well-ordering $\,\preceq\,$.
Assume that for $\;a\in G\;$ 
values $\;\rho(\xi)\in G\,$ are already defined for $\,\xi\prec a\,$.
Put $\;U_a\,:=\,\bigcup_{\xi\prec a}(\rho(\xi)+X_\xi)\,$.
To carry out the inductive step, first of all we make sure 
that the set $\;G\setminus U_a\;$ is not empty.
\smallskip
To verify $\,U_a\not=G\,$
put $\;{\cal T}_a\,:=\,\{\,\rho(\xi)+Y\;|\;x\prec a\,\}\,$.
Then $\,U_a\subset\bigcup{\cal T}_a\subset G\,$.
Let $\,{\cal P}_2\,$ be a partition 
of $\,Y\,$ into sets $\,S\,$ with $\,|S|=2\,$
and put $\,\Sigma\,:=\,\{\,r+s\;|\;\{r,s\}\in{\cal P}_2\,\}\,$.
Since $\,Y\,$ is independent, $\,|\Sigma|=|Y|=\kappa\,$.
We claim that distinct $\,\sigma_1,\sigma_2\in \Sigma\,$ 
never lie in one set $\,T\in{\cal T}_a\,$.
Then in view of $\,|{\cal T}_a|<|\Sigma|\,$
we must have $\;\bigcup{\cal T}_a\not=G\;$ and hence $\,U_a\not=G\,$.
\smallskip
To prove the claim, assume indirectly
$\;\{\sigma_1,\sigma_2\}\,\subset\,\rho(\xi)+Y\;$ 
for $\,\xi\prec a\,$ and distinct $\,\sigma_1,\sigma_2\in \Sigma\,$. 
Then $\;r_1+s_1=\rho(\xi)+y_1\;$
and $\;r_2+s_2=\rho(\xi)+y_2\;$ for distinct $\,y_1,y_2\in Y\,$
and distinct $\,r_1,s_1,r_2,s_2\in Y\,$
and we derive $\;r_1+s_1-r_2-s_2=y_1-y_2\;$
which is obviously impossible since $\,Y\,$ 
is independent.
\smallskip
Now, having made sure that 
$\;G\setminus U_a\;$ is not empty,
we begin the inductive step with the definition 
$\;\theta_a\,:=\,\min_\preceq(G\setminus U_a)\,$.
Distinguishing two cases we put $\,Z_a=X_a\,$ 
if $\,|X_a|=|Y|\,$ and $\,Z_a=Y\,$ if $\,|X_a|<|Y|\,$.
In both cases, $\,X_a\subset Z_a\,$ and $\,|Z_a|=\kappa\,$.
\smallskip
By (5.2) there are exactly $\,\kappa\,$ translates 
$\,t+Z_a\,$ containing the point $\,\theta_a\,$.
Then, trivially, the set 
$\;T\,:=\,\{\,t\in G\;|\;\theta_a\in t\!+\!Z_a\,\}\;$ 
is of size $\,\kappa\,$. 
Since $\;|\{\,\rho(\xi)+X_\xi\;|\;\xi\prec a\,\}|<|T|\,$,
by virtue of (5.1) we can select a set $\,T'\subset T\,$
of size $\,\kappa\,$ such that 
$\;t+Z_a\;$ is disjoint from $\,U_a\,$ 
for every $\,t\in T'\,$. This allows us to conclude the 
inductive step with the definition $\,\rho(a):=\min_\preceq T'\,$.
In doing so, the point $\,\theta_a\,$ lies in 
$\;\rho(a)+X_a\;$ in case that $\,X_a=Z_a\,$.
\smallskip
In this way we obtain pairwise disjoint sets 
$\;\rho(i)+X_i\;(i\in G)\;$ in the group $\,G\,$.
These sets must cover the whole set 
$\,G\,$ because the set 
$\;\{\,i\in G\;|\;\theta_i\in \rho(i)\!+\!X_i\,\}\;$
contains the set $\;J\,:=\,\{\,i\in G\;|\;|X_i|=|Y|\}\;$ 
and hence due to $\,|J|=\kappa\,$ 
the choice function $\,\min_\preceq(\cdot)\,$ 
takes care that no point in $\,G\,$ is accidentally forgotten.
This concludes the proof of Theorem 4.
\medskip\smallskip
{\bf 6.~Bernstein partitions.} In the following $\,{\Bbb P}\,$ 
is an arbitrary uncountable 
Polish space, whence $\,|{\Bbb P}|={\bf c}\,$. 
A {\it Cantor set} is a totally disconnected 
compact metrizable space without isolated points
or, equivalently, a topological space homeomorphic with $\,{\Bbb D}\,$.
Naturally, $\,{\Bbb P}\,$ contains $\,{\bf c}\,$ Cantor sets.
A set $\,B\subset {\Bbb P}\,$ is a {\it Bernstein set}
if and only if no Cantor set is contained in $\,B\,$ or $\,{\Bbb P}\setminus B\,$.
Notice that $\,|B|={\bf c}\,$ for every Bernstein set $\,B\,$.
(One reason for $\,|B|={\bf c}\,$ is (2.1).)
Of course, if $\,{\Bbb P}\,$ is a Banach space then 
translations map Cantor sets onto Cantor sets and hence 
if $\,B\subset {\Bbb P}\,$ is a Bernstein set and $\,\xi\in {\Bbb P}\,$
then $\,\xi+B\,$ is a Bernstein set homeomorphic with $\,B\,$.
Therefore, the following two theorems imply Theorem 3.
\medskip
{\bf Theorem 5.} {\it If $\,\kappa\,$ is a cardinal number with 
$\;2\leq\kappa\leq {\bf c}\;$ and if 
$\,{\Bbb P}\,$ is a separable Banach space
then there exists 
a Bernstein set $\,B\,$ and a set $\,T\subset{\Bbb P}\,$
such that the family $\;{\cal P}\,=\,\{\,t+B\;|\;t\in T\,\}\;$
is a partition of $\,{\Bbb P}\,$ with $\,|{\cal P}|=\kappa\,$.}
\medskip
{\bf Theorem 6.} {\it If $\,\kappa\,$ is a cardinal number
with $\,2\leq \kappa\leq {\bf c}\,$ and
if $\,{\Bbb P}\,$ is an uncountable Polish space then 
$\,{\Bbb P}\,$ can be partitioned into exactly
$\,\kappa\,$ Bernstein sets $\,B\,$ such that 
$\,B\,$ is never homeomorphic to a subspace 
of $\,{\Bbb P}\setminus B\,$.}
\medskip
Since Theorem 6 covers the case that $\,{\Bbb P}\,$ is a Cantor set or 
the Baire space $\,{\Bbb R}\setminus{\Bbb Q}\,$, in comparison with Theorem 3
the following supplement of Theorem 5 is worth mentioning.
\medskip
{\bf Theorem 7.} {\it If $\,{\Bbb P}\,$ is a Cantor set or 
the Baire space $\,{\Bbb R}\setminus{\Bbb Q}\,$ then $\,{\Bbb P}\,$ can be partitioned 
into exactly $\,\kappa\,$ homeomorphic Bernstein sets whenever 
$\,2\leq\kappa\leq{\bf c}\,$.}
\medskip
The case $\,2\leq\kappa<\aleph_0\,$ in Theorem 7 
(and for $\,{\Bbb P}={\Bbb R}\setminus{\Bbb Q}\,$ also the case $\,\kappa=\aleph_0\,$)
is obviously covered by the following noteworthy consequence of Theorem 6. 
\medskip
{\bf Corollary 2.} {\it If for $\,2\leq\lambda\leq\aleph_0\,$
an uncountable Polish space $\,{\Bbb P}\,$ can be partitioned 
into exactly $\,\lambda\,$ homeomorphic open subspaces 
and if either $\,\kappa=\lambda\,$ or $\,2\leq\kappa<\aleph_0=\lambda\,$
then $\,{\Bbb P}\,$ can be partitioned into exactly  $\,\kappa\,$ homeomorphic 
Bernstein sets.}
\medskip
{\it Proof.} In the following, an integer $\,k\geq 0\,$
is divisible by $\,\aleph_0\,$ if and only if $\,k=0\,$. 
Let $\,\{\,U_n\;|\;n\in {\Bbb Z}\,\}\,$ 
be a partition of $\,{\Bbb P}\,$ into $\,\lambda\,$ homeomorphic open subspaces
and let $\,f_n\,$ be a homeomorphism from 
$\,U_0\,$ onto $\,U_n\,$ for each $\,n\in {\Bbb Z}\,$
where $\,U_n=U_m\,$ and $\,f_n=f_m\,$ if and only if 
$\,|n-m|\,$ is divisible by $\,\lambda\,$.
By applying Theorem 6,
let $\,\{\,B_n\;|\;n\in {\Bbb Z}\,\}\,$ be a partition of the Polish space $\,U_0\,$
into $\,\kappa\,$ Bernstein sets where 
$\,B_n=B_m\,$ if and only if 
$\,|n-m|\,$ is divisible by $\,\kappa\,$.
Finally, define 
$\;Y_m\,:=\,\bigcup_{n=-\infty}^\infty f_n(B_{n+m})\;$ for each $\,m\in {\Bbb Z}\,$.
A moment's reflection suffices to see 
that each $\,Y_m\,$ is a Bernstein set in $\,{\Bbb P}\,$.
Clearly, $\;{\cal P}\,=\,\{\,Y_m\;|\;m\in {\Bbb Z}\,\}\;$
is a partition of $\,{\Bbb P}\,$ with $\,|{\cal P}|=\kappa\,$.
Since $\,f_n(B_{n+m})\,$ is always 
an open subset of $\,Y_m\,$, it is evident that $\,Y_m\,$ is homeomorphic 
with $\,Y_0\,$ for each $\,m\in{\Bbb Z}\,$, {\it q.e.d.}
\medskip\smallskip
{\bf 7.~Proof of Theorems 5 and 7.} 
The clue in proving Theorem 5 
is to regard the Banach space
$\,{\Bbb P}\,$ as a vector space over the countable field $\,{\Bbb Q}\,$.
The clue in proving Theorem 7 is to regard 
the Baire space $\,{\Bbb R}\setminus{\Bbb Q}\,$ and the Cantor set $\,{\Bbb D}\,$
as vector spaces over countable fields too.
To justify this viewpoint, let $\,{\Bbb K}\,$ be a countable field
equipped with the {\it discrete} topology.
Naturally, the product space $\,{\Bbb K}^{\aleph_0}\,$
is homeomorphic with the Baire space $\,{\Bbb R}\setminus{\Bbb Q}\,$ 
when $\,{\Bbb K}\,$ is infinite and 
homeomorphic with $\,{\Bbb D}\,$ when $\,{\Bbb K}\,$ is finite.
In all cases, $\,{\Bbb K}^{\aleph_0}\,$ is a vector space 
over $\,{\Bbb K}\,$ of size and dimension $\,{\bf c}\,$.
For a Banach space it is trivial and 
for the vector space $\,{\Bbb K}^{\aleph_0}\,$ it is evident that 
all {\it translations} are {\it continuous}.
In the following we assume that $\,{\Bbb P}\,$ is a Polish space 
and a vector space over a field $\,{\Bbb K}\,$ 
where either $\,|{\Bbb K}|=2\,$ or $\,{\Bbb K}={\Bbb Q}\,$
such that all translations are continuous and 
hence homeomorphisms on $\,{\Bbb P}\,$.
\eject
(If $\,{\Bbb P}\,$ is the Baire space then we do not care 
that $\,{\Bbb Q}\,$ has to be regarded as a
{\it discrete space} since this is done for the only reason 
that $\,{\Bbb P}\,$ can be identified with the product space 
$\,{\Bbb Q}^{\aleph_0}\,$.)
For $\,A\subset{\Bbb P}\,$ let $\,\ell(A)\,$ 
denote the linear subspace of the vector space $\,{\Bbb P}\,$
generated by $\,A\,$. 
It is not difficult to define (by induction)
disjoint Bernstein sets $\;R,S\,\subset\,{\Bbb P}\;$ 
such that $\,R\cup S\,$ is 
{\it linearly independent} in the vector space $\,{\Bbb P}\,$.
(See the remark in Section 8.)
Let $\,\lambda\,$ be a cardinal 
with $\,1\leq\lambda\leq{\bf c}\,$.
Fix a set $\,L\subset S\,$ 
with $\,|L|=\lambda\,$ and let  
$\,H\,$ be a basis of the vector space $\,{\Bbb P}\,$ 
extending the linearly independent set 
$\,R\cup L\,$. Put $\,V=\ell(H\setminus L)\,$ and 
$\,W=\ell(L)\,$. Then the vector space $\,{\Bbb P}\,$ is the direct sum of the 
linear subspaces $\,V\,$ and $\,W\,$. 
\smallskip
The subgroup $\,V\,$ of $\,{\Bbb P}\,$ is a Bernstein set.
Because firstly $\;{\Bbb P}\setminus V\;$ cannot contain 
a Cantor set since $\,R\,$ is a Bernstein set 
and $\;{\Bbb P}\setminus R\;$ contains $\;{\Bbb P}\setminus V\,$.
Therefore, secondly, if some Cantor set $\,C\,$ would lie 
in $\,V\,$ then for every nonzero $\,w\in W\,$ 
the Cantor set $\,w+C\,$ would lie in 
$\,w+V\,$ and hence in $\,{\Bbb P}\setminus V\,$.
Consequently, the quotient group $\;{\Bbb P}/V\,=\,\{\,w+V\,|\;w\in W\,\}\;$
is a partition of $\,{\Bbb P}\,$ into Bernstein sets.
Since the dimension of $\,W\,$ is $\,\lambda\,$,
we have $\;|{\Bbb P}/V|=\lambda\;$
if $\,\lambda\geq\aleph_0\,$. 
So with $\,\kappa=\lambda\,$
the case $\,\aleph_0\leq\kappa\leq{\bf c}\,$ is settled 
in Theorem 5 (with $\,B=V\,$) and in Theorem 7. 
This already concludes the proof of Theorem 7 in view of Corollary 2.
\medskip
In order to conclude the proof of Theorem 5 it remains 
to settle the case $\,2\leq \kappa<\aleph_0\,$ 
for $\,{\Bbb K}={\Bbb Q}\,$. Put $\,\lambda=1\,$ and assume 
$\;\kappa\,\in\,{\Bbb N}\setminus\{1\}\,$. 
Then $\,L\,$ is a singleton $\,\{\xi\}\,$
and $\,W={\Bbb Q}\cdot\xi\,$.
Put $\;M\,:=\,{\Bbb Z}\cap[0,\kappa-1]\;$ 
and $\;Q_\kappa\,:=\,{\Bbb Q}\cap\bigcup_{n\in{\Bbb Z}}[n\kappa,n\kappa+1[\,$.
Then $\;\{\,m+Q_\kappa\;|\;m\in M\,\}\;$  
is a partition of $\,{\Bbb Q}\,$ of size $\,\kappa\,$
and hence $\;\{\,m\cdot\xi+(V+Q_\kappa\cdot\xi)\;|\;m\in M\,\}\;$ 
is a partition of $\,{\Bbb P}\,$ of size $\,\kappa\,$.
Finally, by similar arguments as for proving $\,V\,$ 
to be a Bernstein set,
we realize that $\,V+Q_\kappa\cdot\xi\,$ is a Bernstein set.
Thus the case $\,2\leq \kappa<\aleph_0\,$
in Theorem 5 is settled with $\;B\,=\,V+Q_\kappa\cdot\xi\;$
and $\,T=M\cdot\xi\,$. This concludes the proof  of Theorem 5.
\medskip\smallskip
{\bf 8.~Proof of Theorem 6.} Define a family $\,{\cal F}\,$ 
by $\,f\in{\cal F}\,$ 
if and only if $\,f\,$ is a homeo\-morphism from 
some uncountable G$_\delta$-subset of $\,{\Bbb P}\,$ 
onto some subspace of $\,{\Bbb P}\,$. 
The space $\,{\Bbb P}\,$ is second countable and contains precisely 
$\,{\bf c}\,$ uncountable G$_\delta$-sets.
There are only $\,|{\Bbb P}|^{\aleph_0}={\bf c}\,$ 
mappings from a fixed countable dense 
subset of a Hausdorff space into the set $\,{\Bbb P}\,$.
Consequently, $\,|{\cal F}|={\bf c}\,$
and hence we can  
write $\;{\cal F}\,=\,\{\,f_x\;|\;x\in{\Bbb R}\,\}\;$
without repetitions.
In the following $\,D_x\,$ denotes the domain of $\,f_x\,$
for every $\,x\in{\Bbb R}\,$.
In particular, $\,|D_x|={\bf c}\,$ for every $\,x\in{\Bbb R}\,$.
Clearly, every Cantor set in $\,{\Bbb P}\,$ 
coincides with some G$_\delta$-set $\,D_x\,$. 
(Besides, $\,D_x\,$ always contains a Cantor set.)
Let $\,\preceq\,$ be a minimal well-ordering of the set $\,{\Bbb R}^2\,$.
By induction we define a function $\,g\,$ from
$\,{\Bbb R}^2\,$ into $\,{\Bbb P}\,$ as follows.
Suppose that for $\,(x,y)\in {\Bbb R}^2\,$ 
points $\,g(u,v)\in D_u\,$ 
are already defined for all $\,(u,v)\prec(x,y)\,$.
Then with  
\smallskip
\centerline{$\;E(x,y)\;:=\;
\{\,g(u,v)\;|\;(u,v)\prec(x,y)\,\}\cup
\{\,f_u(g(u,v))\;|\;(u,v)\prec(x,y)\,\}\;$}
\smallskip                                                    
choose any point $\,z\,$ in 
$\;D_x\setminus (E(x,y)\cup f_x^{-1}(E(x,y)))\;$
and define $\,g(x,y):=z\,$.
(This choice is possible since $\,|E(x,y)|<|D_x|\,$ and 
$\,f_x\,$ is injective.) 
In this way the definition of the function $\,g\,$ is finished.
For each $\,t\in {\Bbb R}\,$ define
\smallskip
\centerline{$\;B_t\;:=\;\{\,g(x,t)\;|\;x\in{\Bbb R}\,\}
\cup\{\,f_x(g(x,t))\;|\;x\in{\Bbb R}\,\}\,$.} 
\smallskip
Obviously, if reals $\,s,t\,$ are distinct 
then $\;B_s\cap B_t=\emptyset\,$. 
Furthermore, if $\,t\in {\Bbb R}\,$ then 
$\;D_x\cap B_t\not=\emptyset\;$ for every $\,x\in {\Bbb R}\,$
and hence $\,B_t\,$ meets every Cantor set. 
\smallskip
We claim that for arbitrary $\,t\in {\Bbb R}\,$ 
there does not exist a homeomorphism $\,f\,$ from $\,B_t\,$ 
onto a subspace of $\,{\Bbb P}\setminus B_t\,$.
Suppose indirectly that for some $\,t\in {\Bbb R}\,$ 
such a homeomorphism $\,f\,$ exists. 
Then, by applying the Lavrentieff theorem [2] 4.3.20,
there exists a G$_\delta$-set $\,G\subset{\Bbb P}\,$ with $\,G\supset B_t\,$
and a topological embedding $\,\varphi\,$ from $\,G\,$ into $\,{\Bbb P}\,$ 
such that $\;\varphi(a)=f(a)\;$ for every $\,a\in B_t\,$.
Let $\,x\in{\Bbb R}\,$ be the index where $\,\varphi=f_x\,$
and consider the point $\,a=g(x,t)\,$.
Then $\,a\in B_t\,$ 
and hence $\;f(a)=\varphi(a)=f_x(g(x,t))\in B_t\;$
in contradiction with $\,f(B_t)\cap B_t=\emptyset\,$.
\eject
\smallskip
Finally let $\;2\leq\kappa\leq{\bf c}\;$ and  
fix a set $\,T\subset{\Bbb R}\,$ with $\,0\in T\,$ and $\,|T|=\kappa\,$.
Put $\;U\,:=\,\bigcup\{\,B_t\;|\;0\not=t\in T\,\}\;$
and consider the partition 
\smallskip
\centerline{$\;{\cal P}\,:=\,\{{\Bbb P}\setminus U\}\cup
\{\,B_t\;|\;0\not=t\in T\,\}\;$}
\smallskip
of $\,{\Bbb P}\,$. Trivially, $\,|{\cal P}|=\kappa\,$. 
Since the pairwise disjoint sets 
$\;B_t\;(t\in {\Bbb R})\;$ meet every Cantor set
and since $\;B_0\,\subset\,{\Bbb P}\setminus U\,$,
all elements of $\,{\cal P}\,$ are Bernstein sets.
As already verified, no space $\,B_t\,$ is embeddable into 
the space $\,{\Bbb P}\setminus B_t\,$. 
In particular, $\,B_0\,$ cannot be 
embedded into $\,{\Bbb P}\setminus B_0\,$.
Consequently, $\,{\Bbb P}\setminus U\,$ is not embeddable 
into its complement $\,U\,$. 
Thus $\,{\cal P}\,$ is a partition of $\,{\Bbb P}\,$ 
which proves Theorem 6. 
\medskip
{\it Remark.} The two Bernstein sets $\,R,S\,$ considered
in Section 7
can easily be obtained by simplifying the previous proof 
as follows. Assume that the Polish space $\,{\Bbb P}\,$ is a vector space.
Define $\,\tilde g(x,y)\in{\Bbb P}\,$ inductively 
by picking $\,\tilde g(x,y)\,$ 
from $\;D_x\setminus L(x,y)\;$
where $\,L(x,y)\,$  is the linear subspace of $\,{\Bbb P}\,$ 
generated by all vectors $\,\tilde g(u,v)\,$ with $\,(u,v)\prec (x,y)\,$. 
Put $\;\tilde B_t\,:=\,
\{\,\tilde g(x,t)\;|\;x\in{\Bbb R}\,\}\;$ for $\,t\in{\Bbb R}\,$.
Then $\;\tilde B_t\;(t\in{\Bbb R})\;$ are pairwise disjoint 
Bernstein sets whose union 
is linearly independent.
Thus we can define $\,R:=\tilde B_0\,$ and $\,S:=\tilde B_1\,$.
In a similar way we obtain a Bernstein set $\,B\,$
in the Polish space $\;[0,1]\setminus{\Bbb Q}\;$
where $\,B\subset{\Bbb R}\,$ is linearly independent over $\,{\Bbb Q}\,$.
Then, obviously, $\,{\Bbb Z}+B\,$
is a Bernstein set in the real line. Therefore, 
in view of the proof in Section 7, Theorem 5 remains true 
if the Polish space $\,{\Bbb P}\,$ is the {\it unit circle} $\;x^2+y^2=1\;$
and translations are regarded as additively written {\it rotations}.
\medskip\smallskip
{\bf 9.~More dense partitions.} In the following let 
$\,{\bf P}\,$ denote the class of all Polish spaces without isolated points
which are abelian topological groups 
containing only countably many elements of finite order.
Then all separable Banach spaces and the unit circle 
and the Baire space $\,{\Bbb R}\setminus{\Bbb Q}\,$ lie in $\,{\bf P}\,$.
Furthermore we identify every Cantor set with the topological group
$\,{\Bbb K}^{\aleph_0}\,$ where $\,{\Bbb K}\,$ is the (discrete) field of size $\,3\,$.
(Then the group $\,{\Bbb K}^{\aleph_0}\,$ has no element of order 2.)
Let $\,{\Bbb P}\,$ be a Polish space such that $\,{\Bbb P}\in{\bf P}\,$
or $\,{\Bbb P}\,$ is a Cantor set. (Then $\,|{\Bbb P}|={\bf c}\,$.)
As usual, a set $\,S\subset{\Bbb P}\,$ is $\kappa$-{\it dense} 
if and only if $\;|S\cap U|=\kappa\;$ for every nonempty 
open set $\,U\subset{\Bbb P}\,$.
In view of (2.1) all Bernstein sets in $\,{\Bbb P}\,$ are
${\bf c}$-dense. 
Beside Bernstein sets, other {\it totally imperfect} sets have been
thoroughly investigated. For example the {\it Marczewski sets} 
$\,M\subset{\Bbb P}\,$ defined (in accordance with [4])
such that every Cantor set in $\,{\Bbb P}\,$ 
contains a Cantor set disjoint from $\,M\,$.                      
If $\,\kappa={\bf c}\,$ or $\,\kappa<{\bf c}<2^\kappa\,$ 
then the following theorem allows to 
partition the real line into homeomorphic respectively 
mutually non-homeomorphic $\kappa$-dense sets in various ways.
So it can be accomplished that these $\kappa$-dense sets 
are either meager sets
or Lebesgue null sets or Marczewski sets
including the interesting possibility that (cf.~[4])
the Marczewski sets are {\it Luzin sets} or {\it Sierpi\'nski sets}.
This inclusion is one reason why we consider 
not only ${\bf c}$-density but $\kappa$-density for arbitrary $\,\kappa\leq{\bf c}\,$.
(If a Luzin set $\,L\subset{\Bbb R}\,$
is a Sierpi\'nski set then $\,|L|=\aleph_1\,$, see [4].) 
\medskip
{\bf Theorem 8.} {\it Let $\,{\Bbb P}\,$ be a Polish space in the class 
$\,{\bf P}\cup\{{\Bbb K}^{\aleph_0}\}\,$ 
and let $\,\kappa\geq \aleph_0\,$ and
let $\,{\cal L}\,$ be a $\sigma$-ideal 
in $\,{\Bbb P}\,$ such that $\;t+L\,\in\,{\cal L}\;$
whenever $\,L\in{\cal L}\,$ and $\,t\in{\Bbb P}\,$
and that
every nonempty open set in $\,{\Bbb P}\,$ 
contains a set $\,L\in{\cal L}\,$ with $\,|L|=\kappa\,$.
Then $\,{\Bbb P}\,$ has a partition $\;{\cal P}\subset{\cal L}\;$ 
into $\kappa$-dense sets with $\,|{\cal P}|={\bf c}\,$ 
such that all spaces in $\,{\cal P}\,$ are 
homeomorphic and $\;{\cal P}\,=\,\{\,x+L\;|\;x\in X\,\}\;$
for some $\,L\in{\cal L}\,$ and some $\,X\subset {\Bbb P}\,$.
Furthermore, if $\,2^\kappa>{\bf c}\,$ (particularly if $\,\kappa={\bf c}\,$)
then $\,{\Bbb P}\,$ can also be partitioned into $\,{\bf c}\,$
mutually non-homeomorphic $\kappa$-dense sets lying in $\,{\cal L}\,$.
In case that $\,{\Bbb P}={\Bbb R}\,$ and $\,\kappa>\aleph_0\,$
the real line can be partitioned into $\,{\bf c}\,$
mutually non-homeomorphic dense sets of size $\,\kappa\,$ 
lying in $\,{\cal L}\,$.}
\medskip
{\it Remark.} The condition $\,2^\kappa>{\bf c}\,$ 
rules out $\,\kappa=\aleph_0\,$ in accordance with  
the fact that all countable, 
metrizable spaces without isolated points are homeomorphic with $\,{\Bbb Q}\,$,
see [2] 6.2.A.d.
But even for $\,\kappa>\aleph_0\,$ the condition $\,2^\kappa>{\bf c}\,$
cannot be omitted because 
(see [1]) it is consistent with ZFC that $\,2^{\aleph_1}={\bf c}\,$ and 
{\it all $\aleph_1$-dense subspaces of $\,{\Bbb R}\,$ are homeomorphic.}
\medskip
In order to prove Theorem 8 our first goal 
is to find two disjoint $\kappa$-dense sets in the 
family $\,{\cal L}\,$ whose union is an independent set
in the group $\,{\Bbb P}\,$. Let $\,F({\Bbb P})\,$ be the set of all elements 
of $\,{\Bbb P}\,$ of finite order. (Then $\,F({\Bbb K}^{\aleph_0})={\Bbb K}^{\aleph_0}\,$.)
Put $\,E:=F({\Bbb P})\cup\{0\}\,$ if $\,{\Bbb P}\in{\bf P}\,$ 
and $\,E=\{0\}\,$ if $\,{\Bbb P}={\Bbb K}^{\aleph_0}\,$. Then $\,|E|\leq\aleph_0\,$
in all cases.
For $\,S\subset{\Bbb P}\,$ let $\,[S]\,$ 
denote the subgroup of $\,{\Bbb P}\,$ generated by $\,S\,$.
Let $\,{\cal B}\,$ be a countable basis of the topological space $\,{\Bbb P}\,$.
Choose for each $\,B\in{\cal B}\,$ 
a set $\,L_B\subset B\,$ such that $\,|L_B|={\bf c}\,$ when 
$\,\kappa=\aleph_0\,$ and $\,L_B\in{\cal L}\,$ with $\,|L_B|=\kappa\,$
when $\,\kappa>\aleph_0\,$.
For $\,|A|=\kappa\,$ 
let $\,\preceq\,$ be a
minimal well-ordering of the set 
$\;J\,:=\,A\times{\cal B}\times\{1,2\}\;$ 
of size $\,\kappa\,$.
For $\,j=(a,B,i)\in J\,$ put  $\,\Lambda_j=L_B\,$.
Inductively we define an injective function 
$\,\psi\,$ from $\,J\,$ into $\,{\Bbb P}\,$ such that 
for each $\,y\in J\,$ we pick 
$\,\psi(y)\,$ from $\,\Lambda_y\setminus E\,$ 
satisfying $\;n\cdot\psi(y)\,\not\in\,
[\{\,\psi(x)\;|\;x\prec y\,\}]\setminus\{0\}\;$ 
for every $\,n\in{\Bbb N}\,$. This can be done since 
$\;|\Lambda_y\setminus E|\geq\max\{\aleph_1,\kappa\}>
|[\{\,\psi(x)\;|\;x\prec y\,\}]|\;$  
and $\,|E\!\times\!{\Bbb N}|=\aleph_0\,$.
(Notice that if $\,\{a,b\}\subset{\Bbb P}\in{\bf P}\,$ and $\,n\in{\Bbb N}\,$
then $\,na=nb\,$ implies $\,a\!-\!b\in E\,$ 
and that all subgroups of $\,{\Bbb K}^{\aleph_0}\,$ are linear subspaces.)
It is evident that 
$\;L_i\,:=\,\psi(A\times{\cal B}\times\{i\})\;$ 
is a $\kappa$-dense subset of $\,{\Bbb P}\,$ for $\,i\in\{1,2\}\,$
and $\;L:=\psi(J)=L_1\cup L_2\;$ is independent in the group $\,{\Bbb P}\,$.
By definition, $\,L\cap E=\emptyset\,$.
If $\,\kappa=\aleph_0\,$ then $\,L\,$
is a countable set and hence $\,L\in{\cal L}\,$. 
If $\,\kappa>\aleph_0\,$ then $\,L\in{\cal L}\,$
since $\;L\,\subset\,\bigcup\,\{\,L_B\;|\;B\in{\cal B}\,\}\;$
and $\,L_B\in{\cal L}\,$ for all $\,B\,$ in the countable 
family $\,{\cal B}\,$. 
\smallskip
Now put $\;{\cal X}\,:=\,\{\,L_1\cup S\;|\;S\subset L_2\,\}\,$.
The family $\,{\cal X}\,$ 
consists of $\kappa$-dense sets from the family $\,{\cal L}\,$
lying in $\,L\,$. Trivially, $\,|{\cal X}|=2^\kappa\,$.
If $\,2^\kappa>{\bf c}\,$ then  
there must exist a family $\,{\cal Y}\subset {\cal X}\,$
with $\,|{\cal Y}|=|{\cal X}|\,$ such that 
the spaces in $\,{\cal Y}\,$ are mutually non-homeomorphic.
(Because if $\,H\,$ is 
a second countable Hausdorff space then 
there are at most $\,{\bf c}\,$ continuous functions from a 
subspace of $\,H\,$ into $\,H\,$.) 
Consequently, we have constructed one
$\kappa$-dense set $\,L\in{\cal L}\,$
independent in the group $\,{\Bbb P}\,$ with $\,a+a\not=0\,$ for 
every $\,a\in L\,$.
If $\,2^\kappa>{\bf c}\,$ then $\,L\,$ has 
at least $\,{\bf c}\,$ mutually non-homeomorphic 
$\kappa$-dense subsets. 
Furthermore we claim that if $\,{\Bbb P}={\Bbb R}\,$ and 
$\,\kappa>\aleph_0\,$ then $\,L\,$  
has $\,{\bf c}\,$ mutually non-homeomorphic dense subsets
of size $\,\kappa\,$. 
\smallskip
Indeed, we can create sets $\,D\subset L\,$ of size $\,\kappa\,$
such that 
for certain $\,u<v\,$ the set $\,D\cap[u,v]\,$ is {\it black}
respectively {\it grey}, meaning that 
$\,|D\cap[x,y]|=\kappa\,$, respectively $\,|D\cap[x,y]|=\aleph_0\,$,
whenever $\,u\leq x<y\leq v\,$.
There is enough freedom for coloring the real line 
in grey and black sections in order to 
create $\,{\bf c}\,$ mutually 
non-homeomorphic sets $\,D\,$ of that kind. Actually, this 
can be accomplished in the same 
way as in the proof of [5] Theorem~6.
(To give a hint here, if $\,u<v<w\,$
and $\,D\cap[u,v]\,$ is grey and $\,D\cap[v,w]\,$ is black
then the point $\,v\,$ is the unique point $\,a\,$ in 
$\;D\cap[u,w]\;$ satisfying the topological property 
that every neighborhood of $\,a\,$ 
contains both some uncountable open set and some
countable open set.)
Considering all these possibilities concerning 
subsets of the independent $\kappa$-dense set $\,L\in{\cal L}\,$,
Theorem 8 is settled by applying Theorem 4 
with $\,G={\Bbb P}\,$ and $\,Y=L\,$.           
\bigskip\bigskip
{\bf References} 
\medskip\smallskip
\vbox{\eightpoint [1] Baumgartner, J.: 
{\it All $\aleph_1$-dense sets of reals can be isomorphic.} 
Fund.~Math.~{\bf 79} (1973) 101-106.
\smallskip
[2] Engelking, R.: {\it General 
Topology, revised and completed edition.} Heldermann 1989. 
\smallskip 
[3] Jech, T.: {\it Set Theory.} 3rd ed. Springer 2002. 
\smallskip
[4] Kharazishvili, A.B.: 
{\it Strange Functions in Real Analysis.} Chapman and Hall 2006. 
\smallskip
[5] Kuba, G.: {\it Counting linearly ordered spaces.}
Colloq.~Math.~{\bf 135} (2014), 1-14.
\smallskip 
[6] Kuba, G.: 
{\it Counting ultrametric spaces.} Colloq.~Math.~{\bf 152} (2018), 217-234. 
\smallskip
[7] Rotman, J.: {\it An Introduction to the Theory of Groups.}  
Springer 1995. 
\smallskip
[8] Wagon, S.: {\it The Banach-Tarski Paradox.} Cambridge 
University Press 1986. 
\bigskip
Gerald Kuba \smallskip Institute of 
Mathematics, 

University of Natural Resources and Life Sciences, 1180 Wien, Austria. 
\smallskip {\sl E-mail:} {\tt gerald.kuba(at)boku.ac.at}} 
\end